\documentclass[reqno]{amsart}
\usepackage{amssymb,latexsym}

\newcommand{\bfi}{\bfseries\itshape}

\newtheorem{thm}{Theorem}[section]

\newtheorem{cor}{Corollary}[section]

\newtheorem{defn}{Definition}[section]

\newtheorem{rem}{Remark}[section]

\makeatletter
\@addtoreset{figure}{section}
\def\thefigure{\thesection.\@arabic\c@figure}
\def\fps@figure{h, t}
\@addtoreset{table}{bsection}
\def\thetable{\thesection.\@arabic\c@table}
\def\fps@table{h, t}
\@addtoreset{equation}{section}

\makeatother

\allowdisplaybreaks
\def\intprod{\mathbin{\hbox to 6pt{%
          \vrule height0.4pt width5pt depth0pt
          \kern-.4pt
          \vrule height6pt width0.4pt depth0pt\hss}}}

\begin{document}

\title
[Reduction of Discrete Lagrangians on Lie groups]
{Symmetry Reduction of Discrete Lagrangian Mechanics on Lie groups}

\author[J.E. Marsden]{Jerrold E. Marsden}
\address[Marsden]
{Control and Dynamical Systems\\California Institute of Technology,
  107-81\\ Pasadena, CA 91125}
\email{marsden@cds.caltech.edu}

\author[S. Pekarsky]{Sergey Pekarsky}
\address[Pekarsky]
{Control and Dynamical Systems\\California Institute of Technology,
  107-81\\ Pasadena, CA 91125}
\email{sergey@cds.caltech.edu}

\author[S. Shkoller]{Steve Shkoller}
\address[Shkoller]{Department of Mathematics\\
University of California\\
Davis, CA 95616}
\email{shkoller@math.ucdavis.edu}

\subjclass{Primary 70H35,70E15; Secondary 58F}

\date{February 1999; current version \today}
\keywords{Euler-Poincar\'{e}, symplectic, Poisson}

\begin{abstract}
For a discrete mechanical system on a Lie group $G$
determined by a (reduced) Lagrangian $\ell$ we define a
Poisson structure via the pull-back of the Lie-Poisson
structure on the dual of the Lie algebra
${\mathfrak g}^*$ by the corresponding Legendre
transform. The main result shown in this paper is that this
structure coincides with the reduction under the symmetry group $G$
of the canonical discrete Lagrange $2$-form
$\omega_\mathbb{L}$ on
$G \times G$. Its symplectic leaves then become dynamically
invariant manifolds for the reduced discrete system. Links between our
approach and that of groupoids and algebroids as well as the reduced
Hamilton-Jacobi equation are made. The rigid body is discussed as an
example.
\end{abstract}

\maketitle

\tableofcontents

%%%%%%%%%%%%%%%%%%%%%%%%%%%%%%%%%%%%%%%%%%%%%%%%
\section{Introduction}
%%%%%%%%%%%%%%%%%%%%%%%%%%%%%%%%%%%%%%%%%%%%%%%%

\noindent {\bf Background.} This paper continues our
development of discrete Lagrangian mechanics on a Lie group
introduced in Marsden, Pekarsky, and Shkoller \cite{MPeS}. In
our earlier paper, using the context of the Veselov method
for discrete mechanics, discrete analogues of
Euler-Poincar\'{e} and Lie-Poisson reduction theory (see, for example,
Marsden and Ratiu \cite{MaRa1999}) were
developed for systems on finite dimensional Lie groups $G$
with Lagrangians $L:TG \rightarrow  \mathbb{R}$ that are
$G$-invariant. The resulting discrete equations provide
``reduced'' numerical algorithms which manifestly preserve
the symplectic  structure. The manifold $G \times G$ is used
as the discrete approximation of $TG$, and a discrete
Lagrangian $\mathbb{L}:G \times G \rightarrow
\mathbb{R}$  is constructed from a given Lagrangian $L$
in such a way that the
$G$-invariance property is preserved.  Reduction by $G$
results in a new ``variational'' principle for the reduced
Lagrangian $\ell: G \cong (G \times G) /G \rightarrow \mathbb{R}$,
which then
determines the discrete  Euler-Poincar\'{e} (DEP) equations.
Reconstruction of these equations is consistent with the
usual Veselov discrete Euler-Lagrange equations developed in
\cite{WM,MPS},  which are naturally symplectic-momentum
algorithms. Furthermore, the solution of the DEP algorithm
leads directly to a  discrete Lie-Poisson (DLP) algorithm.
For example, when $G=\text{SO}(n)$, the DEP and DLP
algorithms for a particular choice of  the discrete
Lagrangian $\mathbb{L}$ are equivalent to the Moser-Veselov
\cite{MoV}  scheme for the generalized rigid body.
\medskip

\noindent {\bf Main Results of this Paper.} We show that when
a discrete Lagrangian $\mathbb{L}: G \times G \rightarrow
\mathbb{R}$ is $G$-invariant, a Poisson structure on (a
subset) of one copy of the Lie group $G$ can be defined which
governs the corresponding discrete reduced dynamics.  The
symplectic leaves of this structure become dynamically
invariant manifolds which are manifestly preserved under the
structure preserving discrete Euler-Poincar\'{e} algorithm
(see Section 2.1).

Moreover, starting with a discrete Euler-Poincar\'{e}
system on $G$ one can readily recover, by
means of the Legendre transformation, the corresponding
Lie-Poisson  Hamilton-Jacobi system on ${\mathfrak g}^*$
analyzed by Ge and Marsden \cite{GM}; the relationship
between the discrete Euler-Lagrange and discrete
Euler-Poincar\'{e} equations and the Lie-Poisson
Hamilton-Jacobi equations was  examined from a different
point of view in our companion paper \cite{MPeS}.

We also apply Weinstein's results on
Lagrangian mechanics on groupoids  and algebroids \cite{W96}
to the setting of regular Lie groups.  The groupoid-algebroid
setting reveals new and  interesting connections between
discrete and continuous dynamics.

%%%%%%%%%%%%%%%%%%%%%%%%%%%%%%%%%%%%%%%%%%%%%%%%
\section{Discrete reduction}
%%%%%%%%%%%%%%%%%%%%%%%%%%%%%%%%%%%%%%%%%%%%%%%%

In this section we review the discrete
Euler-Poincar\'{e} reduction of a Lagrangian system on $G\times G$
considered in detail in \cite{MPeS}.
We approximate $TG$ by $G \times G$ and form a discrete Lagrangian
$\mathbb{L}: G \times G \rightarrow \mathbb{R}$ from the original
Lagrangian $L:TG \rightarrow \mathbb{R}$  by
$$
\mathbb{L}(g_k, g_{k+1}) = L (\kappa (g_k, g_{k+1}),
\chi (g_k, g_{k+1}) ),
$$
where $\kappa$ and $\chi$ are functions of $(g_k, g_{k+1})$
which approximate the current configuration $g(t) \in G$ and the corresponding
velocity $\dot{g}(t) \in T_g G$.
We choose discretization schemes  for which the
discrete Lagrangian $\mathbb{L}$ inherits the symmetries of the original
Lagrangian $L$:  $\mathbb{L}$ is $G$-invariant on $G\times G$
whenever $L$ is $G$-invariant on $TG$.
In particular, the induced right (left) lifted action of $G$ onto $TG$
corresponds to the diagonal right (left) action of $G$ on
$G \times G$.

Having specified the discrete Lagrangian, we form the \emph{action sum}
$$
\mathbb{S} = \sum_{k=0}^{N-1} \mathbb{L}(g_k, g_{k+1}),
$$
which approximates the action integral $S = \int L dt$,
and obtain the discrete Euler-Lagrange (DEL) equations
\begin{equation}
   \label{DEL}
   D_2 \mathbb{L}(g_{k-1}, g_k) + D_1 \mathbb{L}(g_k, g_{k+1}) = 0,
\end{equation}
as well as the discrete symplectic  form $\omega_{\mathbb{L}}$,
given in coordinates on $G \times G$ by
%\begin{equation}\label{dsymp}
$$
\omega_\mathbb{L}= \frac{\partial^2 \mathbb L}{\partial g^i_k
\partial g^j_{k+1}} dg_k^i \wedge dg_{k+1}^j.
$$
%\end{equation}
In (\ref{DEL}), $D_1$ and $D_2$ denote derivatives with respect
to the first and second argument, respectively.
The algorithm (\ref{DEL}) as well as $\omega_\mathbb{L}$ are obtained
by extremizing the action sum $\mathbb{S} : G^{N+1} \rightarrow \mathbb{R}$
with
arbitrary variations.
Using this variational point of view, it is known  that the
flow $\mathbb{F}_t$ of the DEL equations preserves this discrete
symplectic structure. This result was obtained  using a
discrete Legendre transform and a direct computation in
\cite{V88,V91,WM} and  a proof using the variational structure directly
was given in \cite{MPS}.

\begin{rem}
\label{rem2.1}
{\rm We remark that the discrete symplectic structure
$\omega_\mathbb{L}$  is not globally defined, but rather need only be
nondegenerate in a neighborhood 
of the diagonal $\Delta$ in $G \times G$, i.e., whenever
$g_k$ and $g_{k+1}$ are nearby.  Section 3 of
\cite{MPS} shows that $\omega_\mathbb{L}$ arises from the boundary
terms of the discrete action sum restricted to the space of solutions
of the discrete Euler-Lagrange equations; an implicit function theorem
argument relying on the regularity of the discrete Lagrangian $\mathbb{L}$
is required in order to obtain solutions to the discrete Euler-Lagrange
equations,
and this regularity need only hold in a neighborhood of the
diagonal $\Delta \subset G\times G$.}
\end{rem}

%%%%%%%%%%%%%%%%%%%%%%%%%%%%%%%%%%%%%%%%%%%%%%%%
\subsection{The discrete Euler-Poincar\'{e} algorithm}
%%%%%%%%%%%%%%%%%%%%%%%%%%%%%%%%%%%%%%%%%%%%%%%%
The discrete reduction of a right-invariant system proceeds
as follows (see \cite{MPeS} for details). The case of left invariant
systems is similar. Of course, some
systems such as the rigid body are left invariant.

The induced group action on $G \times G$ by
an element $\bar{g} \in G $ is simply right multiplication in
each component:
$$
\bar{g} : (g_k, g_{k+1}) \mapsto (g_k \bar{g}, g_{k+1} \bar{g}),
$$ for all $g_k, g_{k+1} \in G.$

The quotient map is given by
\begin{equation}\label{projection}
\pi_d : G \times G \rightarrow (G \times G)/G \cong G,\
\quad (g_k, g_{k+1}) \mapsto g_k g_{k+1}^{-1}.
\end{equation} One may alternatively use
$g_{k+1} g_k^{-1}$ instead of
$g_k g_{k+1}^{-1}$ as the quotient map; the projection map
(\ref{projection}) defines the
{\bfi reduced discrete Lagrangian}
$\ell : G \rightarrow \mathbb{R}$ for any $G$-invariant
$\mathbb{L}$ by $\ell \circ \pi_d = \mathbb{L}$, so that
$$ \ell(g_k g_{k+1}^{-1}) = \mathbb{L} (g_k, g_{k+1}), $$ and
the {\bfi reduced action sum} is given by
$$ s = \sum_{k=0}^{N-1} \ell (f_{k k+1}), $$ where $f_{k k+1}
\equiv g_k g_{k+1}^{-1}$ denotes points in the quotient space.
A reduction of the DEL equations results in the  {\bfi
discrete Euler-Poincar\'{e}} (DEP) equations
\begin{equation} \label{DEP}
R^*_{f_{k k+1}}\ell' (f_{k k+1}) -
L^*_{f_{k-1 k}} \ell' (f_{k-1 k}) = 0
\end{equation}
for $k=1,...,N-1$,
where $R^*_f$ and $L^*_f$, for $f \in G $ are the right and left
pull-backs by $f$, respectively, defined as follows: for $\alpha _g \in T
^\ast _g G $, $R^*_f \alpha_g \in \mathfrak{g}^\ast$ is given by
$\langle R^*_f \alpha_g, \xi \rangle =
\left\langle \alpha _g, T R _f \cdot \xi \right\rangle$ 
for any $\xi \in {\mathfrak g}$, where
$T R _f $ is the tangent map of the right translation map $R _f: G
\rightarrow G $; $h \rightarrow h f $, 
with a similar definition for $L^*_f$. Also,
$\ell^\prime: G \rightarrow T^*G$ is the differential of $\ell$ defined
as  follows.
Let $g^\epsilon$ be a smooth curve in $G$ such that $g^0 =g$ and
$(d/d\epsilon)|_{\epsilon = 0} g^\epsilon = v$.  Then
$$
\ell^\prime (g)\cdot v= (d/d\epsilon)|_{\epsilon = 0} \ell(g^\epsilon ).
$$
For the other choice of the quotient in (\ref{projection}) given
by $h_{k+1 k} \equiv g_{k+1} g_k^{-1}$, the DEP equations are
\begin{equation} \label{DEP2}
L^*_{h_{k+1 k}}\ell' (h_{k+1 k}) -
R^*_{h_{k k-1}} \ell' (h_{k k-1}) = 0
\end{equation}

\begin{rem}
{\rm In the case that $\mathbb{L}$ is left invariant, the
discrete Euler-Poincar\'{e} equations take the form
\begin{equation}\label{DEPleft}
L^*_{f_{k+1 k}}\ell' (f_{k+1 k}) -
R^*_{f_{k k-1}} \ell' (f_{k k-1}) = 0
\end{equation}
where  $f_{k+1k} \equiv g^{-1}_{k+1}g_k$ is in the left quotient
$(G\times G)/G$.

Notice that equations
(\ref{DEP2}) and (\ref{DEPleft}) are formally the same.}
\end{rem}

We may associate to any $C^1$ function $F$ defined on a neighborhood
$\mathcal{V}$ of
$\Delta \subset G \times G$ its Hamiltonian
vector field $X_F$ on $ \mathcal{V} \supset \Delta$ satisfying
$X_F \intprod \omega_\mathbb{L} = dF$, where $d F$, the differential of
$F$, is a one-form. The symplectic structure $\omega_\mathbb{L}$ naturally
defines a Poisson structure  on a neighborhood $\mathcal{V}$ of  $\Delta$
(which we shall denote
$\{ \cdot, \cdot \}_{G\times G}$) by the usual relation
$$
\{F,H\}_{G \times G} = \omega_\mathbb{L}(X_F,X_H).
$$

Theorem $2.2$ of \cite{MPeS} states that if the action of $G$ on $G\times
G$ is proper,  the algorithm on $G$ defined by the discrete
Euler-Poincar\'{e} equations (\ref{DEP}) preserves the induced Poisson
structure $\{ \cdot, \cdot \}_G$ on
$\mathcal{U} \subset G$ given by
\begin{equation}\label{p2}
\{f,h\}_G \circ \pi_d = \{ f \circ \pi_d, h \circ \pi_d \}_{G \times G}
\end{equation}
for any $C^1$ functions $f,h$ on $\mathcal{U}$, where
$\mathcal{U} = \pi_d (\mathcal{V})$.

Using the  definition $f_{k k+1}=g_kg_{k+1}^{-1}$, the
DEL algorithm can be  reconstructed from the DEP algorithm by
\begin{equation}
(g_{k-1}, g_k) \mapsto (g_k, g_{k+1}) =
(f^{-1}_{k-1 k} \cdot g_{k-1}, f^{-1}_{k k+1} \cdot g_k),
\label{rec_DEL}
\end{equation}
where $f_{k k+1  }$ is the solution of (\ref{DEP}).
Indeed, $f^{-1}_{k k+1} \cdot g_k$ is precisely $g_{k+1}$.
Similarly one shows that in the  case of a left $G$ action,
the reconstruction of the DEP equations (\ref{DEPleft})
is given by
%\begin{equation}\label{ref_DEL_left}
$$
(g_{k-1}, g_k) \mapsto (g_k, g_{k+1}) =
(g_{k-1} \cdot f^{-1}_{k k-1}, g_k \cdot f^{-1}_{k+1 k}).
$$
%\end{equation}

%%%%%%%%%%%%%%%%%%%%%%%%%%%%%%%%%%%%%%%%%%%%%%%%
\subsection{The discrete Lie-Poisson algorithm}
%%%%%%%%%%%%%%%%%%%%%%%%%%%%%%%%%%%%%%%%%%%%%%%%

In addition to reconstructing the dynamics on  $G\times G$, one may use
the coadjoint action to form a {\bfi discrete Lie-Poisson} algorithm
approximating the dynamics on
$\mathfrak{g}^\ast$
\cite{MPeS}
\begin{equation}
   \label{DLP}
   \mu_{k+1} = \operatorname{Ad}^\ast_{f_{k k+1}} \cdot \mu_k,
\end{equation}
where $\mu_k := \operatorname{Ad}^\ast_{g^{-1}_k} \mu_0$
is an element of the dual of the Lie algebra,
$\mu_0$ is the constant of motion (the momentum map value), and
the sequence $\{f_{kk+1}\}$ is provided by the DEP algorithm on $G$.

The corresponding discrete Lie-Poisson equation for the left invariant
system is given by
\begin{equation}\label{DLP_left}
\Pi_{k+1} = \operatorname{Ad}^\ast_{f^{-1}_{k+1 k}} \cdot \Pi_k,
\end{equation}
where $\Pi_k := \operatorname{Ad}^\ast_{g_k} \pi_0$ and
$\pi_0$ is the constant momentum map value.
Henceforth, we shall use the notation
$\mu \in  {\mathfrak{g}}^\ast$ for the \emph{right} invariant system
and $\Pi \in  {\mathfrak{g}}^\ast$ for the \emph{left}.

%%%%%%%%%%%%%%%%%%%%%%%%%%%%%%%%%%%%%%%%%%%%%%%%%%%%%%%%%
\section{Poisson structure and invariant manifolds on Lie groups}
%%%%%%%%%%%%%%%%%%%%%%%%%%%%%%%%%%%%%%%%%%%%%%%%%%%%%%%%%

Discretization of an Euler-Poincar\'{e} system on $T G$ results in a
system on $G \times G$ defined by a Lagrangian $\mathbb{L}$.
If it is regular, the Legendre transformation (in the sense of Veselov) $F
\mathbb{L}$ define a symplectic form (and, hence, a Poisson structure)
on $\mathcal{V} \subset G \times G$ via the pull-back of the
canonical form from $T^*G$.
Then, general Poisson reduction applied to these discrete
settings defines a  Poisson structure on the reduced
space $\mathcal{U} = \pi_d (\mathcal{V}) \subset G$.
This approach was adopted in Theorem 2.2 of \cite{MPeS}.

Alternatively, without appealing to the reduction procedure,
a Poisson structure on a Lie group can be defined using
ideas of Weinstein \cite{W96} on Lagrangian mechanics on
groupoids and their algebroids. The key idea can be summarized
in the following statements. A smooth function on a groupoid
defines a natural (Legendre type) transformation between
the groupoid and the dual of its algebroid. This transformation
can be used to pull back a canonical Poisson structure from the dual
of the algebroid, provided the regularity conditions are satisfied.

The ideas outlined in this section can be easily expressed using the
groupoid-algebroid formalism. Such a formalism is suited to the discrete
gauge field theory generalization as well as to
discrete semidirect product theory;
nevertheless, the theory of groupoids and algebroids
is not essential for the derivations, but rather
contributes nicely to the elegance of the exposition.

%%%%%%%%%%%%%%%%%%%%%%%%%%%%%%%%%%%%%%%%%%%%%%%%%%%%%%%%%
\subsection{Dynamics on groupoids and algebroids}
%%%%%%%%%%%%%%%%%%%%%%%%%%%%%%%%%%%%%%%%%%%%%%%%%%%%%%%%%

In this subsection, we show that our discrete reduction methodology is
consistent with Weinstein's groupoid-algebroid construction;  the contents
of this subsection are not essential for the remainder of the paper.

We briefly summarize results from Weinstein \cite{W96} and
refer the reader to the original paper for details of proofs
and definitions.
Let $\Gamma$ be a groupoid over a set $M$, with
$\alpha, \beta : \Gamma \rightarrow M$ being its source and target
maps, with a multiplication map $m : \Gamma_2 \rightarrow \Gamma$,
where $\Gamma_2 \equiv \{ (g,h) \in \Gamma \times \Gamma \ | \
\beta(g) = \alpha(h) \}$.
Denote its corresponding algebroid by ${\mathcal A}$.

The Lie groupoids relevant to our exposition are the Cartesian
product $G \times G$ of a Lie group $G$, with multiplication
$(g,h)(h,k) = (g,k)$, and the group $G$ itself.
The corresponding algebroids are the tangent bundle $TG$
and the Lie algebra ${\mathfrak g}$, respectively.
The dual bundle to a Lie algebroid carries a natural Poisson
structure. This is the Poisson bracket associated to the
canonical symplectic form on $T^*G$ and the Lie-Poisson
structure on ${\mathfrak g}^*$, respectively.

Lagrangian mechanics on a groupoid $\Gamma$ is defined as follows.
Let ${\mathcal L}$ be a smooth, real-valued function on $\Gamma$,
${\mathcal L}_2$ the restriction to $\Gamma_2$ of the function
$(g,h) \mapsto {\mathcal L}(g) + {\mathcal L}(h)$.

\begin{defn}
\label{def3.1}
Let $\Sigma_{\mathcal L} \subset \Gamma_2$ be the set of
critical points of ${\mathcal L}_2$ along the fibers of the
multiplication map $m$;
i.e. the points in $\Sigma_{\mathcal L}$
are stationary points of the function ${\mathcal L}(g)+{\mathcal L}(h)$
when $g$ and $h$ are restricted to admissible pairs with the constraint
that the product $g h$ is fixed \cite{W96}.
\end{defn}
A {\bfi solution of the Lagrange equations} for the
Lagrangian function
${\mathcal L}$ is a sequence $\dots, g_{-2}, g_{-1}, g_0, g_1, g_2,
\dots$ of elements of $\Gamma$, defined on some ``interval''
in ${\mathbb Z}$, such that $(g_j, g_{j+1}) \in \Sigma_{\mathcal L}$
for each $j$.

The Hamiltonian formalism for discrete Lagrangian systems is based
on the fact that each Lagrangian submanifold of a symplectic groupoid
determines a Poisson automorphism on the base Poisson manifold.
Recall that the cotangent bundle $T^* \Gamma$ is, in addition
to being a symplectic manifold, a groupoid itself, the base
being ${\mathcal A}^*$; notice that both manifolds are naturally
Poisson. The source and target mappings
$\tilde{\alpha}, \tilde{\beta} : T^*\Gamma \rightarrow {\mathcal A}^*$
are induced by $\alpha$ and  $\beta$.
\begin{defn}
\label{def3.2}
Given  any smooth function ${\mathcal L}$ on $\Gamma$, a
Poisson map
$\Lambda_{\mathcal L}$ from ${\mathcal A}^*$ to itself,
which may be said to be generated by ${\mathcal L}$
is defined by the Lagrangian submanifold $d {\mathcal L}(\Gamma)$
(under a suitable hypothesis of nondegeneracy) \cite{W96}.
\end{defn}
The appropriate ``Legendre transformation'' $F {\mathcal L}$ in the
groupoid context is given by $\tilde\alpha \circ d {\mathcal L} :
\Gamma \rightarrow {\mathcal A}^*$ or $\tilde\beta \circ d {\mathcal L} :
\Gamma \rightarrow {\mathcal A}^*$, depending on whether we
consider  right or left invariance (through the definition
of maps $\tilde\alpha$ and $\tilde\beta$).
The transformation $F {\mathcal L}$ relates the mapping on $\Gamma$
defined by $\Sigma_{\mathcal L}$ with the mapping
$\Lambda_{\mathcal L}$ on ${\mathcal A}^*$.
$F {\mathcal L}$ also pulls back the Poisson structure from
${\mathcal A}^*$ to $\Gamma$, which, in general, is defined
only \emph{locally} on some neighborhood $\mathcal{U} \subset \Gamma$.
In the context of a Lie group, this means that any regular function
$\ell : G \rightarrow \mathbb{R}$ defines a Poisson structure
on $\mathcal{U}$. We shall address this issue in the next subsections.
The reader is referred to \cite{W96} for an application
of the above ideas to the groupoid $M \times M$ when the manifold $M$
does not necessarily have group structure.

%%%%%%%%%%%%%%%%%%%%%%%%%%%%%%%%%%%%%%%%%%%%%%%%%%%%%%%%%
\subsection{DEP equations as generators of Lie-Poisson Hamilton-Jacobi
   equations}
%%%%%%%%%%%%%%%%%%%%%%%%%%%%%%%%%%%%%%%%%%%%%%%%%%%%%%%%%

A Lie group $G$ is the simplest example of a groupoid with the base
being just a point. Its algebroid is the corresponding Lie
algebra ${\mathfrak g}$, with the dual being ${\mathfrak g}^*$.
Consider left invariance and let a general function ${\mathcal L}$
on the group be specified by the discrete reduced Lagrangian
$\ell : G \rightarrow \mathbb{R}$.
Then, the Legendre transform defined above is given by
%\begin{equation} \label{legendre}
$$
F \ell = L^*_g \circ d \ell \ : G \rightarrow {\mathfrak g}^*,
$$
%\end{equation}
where $d \ell : G \rightarrow T^* G$.
Using these transformations we define
$$
\Pi_{k-1} \equiv F\ell (f_{k k-1})=L^*_{f_{k k-1}} \circ d \ell (f_{k k-1}).
%\Pi_{k+1} \equiv F\ell (f_{k+1 k})=L^*_{f_{k+1 k}} \circ d \ell ({k+1 k}).
$$
Recall the DEP equation (\ref{DEPleft}) for left-invariant systems :
$$
L^*_{f_{k+1 k}} d \ell (f_{k+1 k}) -
R^*_{f_{k k-1}} d \ell (f_{k k-1}) = 0,
$$
where we have identified the notations $\ell'$ and $d \ell$.
The latter equation can be rewritten as a system
\begin{equation}
\label{DEPHJ}
\left\{ \begin{array}{l}
\Pi_{k} = L^*_f \circ d \ell (f), \\
\Pi_{k+1} = R^*_f \circ d \ell (f),
\end{array} \right.
\end{equation}
where the first equation is to be solved for $f$ (which stands
for $f_{k+1 k}$) which then is substituted into the second
equation to compute $\Pi_{k+1}$.

This system is precisely the Lie-Poisson Hamilton-Jacobi system
described in \cite{GM} with the reduced discrete Lagrangian $\ell$
playing the role of the generating function. This means that there is
no need to find an approximate solution of the reduced Hamilton-Jacobi
equation \cite{GM}.
Notice also that the DLP equation (\ref{DLP_left}) is a direct
consequence of the system (\ref{DEPHJ}):
$$
\Pi_{k+1} = \operatorname{Ad}^\ast_{f^{-1}_{k+1 k}} \cdot \Pi_k.
$$

The following diagrams relate the dynamics on $G$ and on ${\mathfrak g}^*$:
\begin{equation}
\label{diag_dyn}
\begin{array}{ccc}
G & \stackrel{\Sigma_{\ell}}{\longrightarrow} & G \\
\Big\downarrow\vcenter{%
   \rlap{$\scriptstyle{\rm }\,F \ell$}} &
         & \Big\downarrow\vcenter{%
            \rlap{$\scriptstyle{\rm }\,F \ell$}}\\
{\mathfrak g}^* & \stackrel{\Lambda_{\ell}}{\longrightarrow} & {\mathfrak g}^*
\end{array}
\qquad \qquad
\begin{array}{ccc}
f_{k k-1} & \stackrel{\Sigma_{\ell}}{\longrightarrow} & f_{k+1 k} \\
\Big\downarrow\vcenter{%
   \rlap{$\scriptstyle{\rm }\,F \ell$}} &
         & \Big\downarrow\vcenter{%
            \rlap{$\scriptstyle{\rm }\,F \ell$}}\\
\Pi_{k-1} & \stackrel{\Lambda_{\ell}}{\longrightarrow} & \Pi_{k},
\end{array}
\end{equation}
where $\Sigma_\ell$ and $\Lambda_\ell$ are given in
Definitions \ref{def3.1} and \ref{def3.2}.

%%%%%%%%%%%%%%%%%%%%%%%%%%%%%%%%%%%%%%%%%%%%%%%%%%%%%%%%%
\subsection{Some Advantages of Structure-preserving
Integrators}
%%%%%%%%%%%%%%%%%%%%%%%%%%%%%%%%%%%%%%%%%%%%%%%%%%%%%%%%%

As we mentioned above, the ``Legendre transform'' $F \ell$
allows us to put a Poisson structure on the Lie group $G$,
which, of course,  depends on the discrete Lagrangian
${\mathbb L}$ on $G\times G$,
and hence on the original Lagrangian $L$ on $TG$
(if we consider this from the discrete reduction point of view).
It follows that the reduction of the discrete Euler-Lagrange dynamics
on $G \times G$ is necessarily restricted
to the symplectic leaves of this Poisson structure, so that these leaves
are invariant manifolds,  and correspond (under $F \ell$)
to the symplectic leaves (coadjoint orbits) of the continuous reduced
system on ${\mathfrak g}^*$.

These ideas are the content of the following theorems.

\begin{thm}
Let $L$ be a right invariant Lagrangian on
$TG$ and let $\mathbb{L}$ be the Lagrangian of the
corresponding discrete system on $\mathcal{V} \subset G \times G$.
Assume that $\mathbb{L}$ is regular, in the sense that
the Legendre transformation $F \mathbb{L} : \mathcal{V} \rightarrow
F \mathbb{L} (\mathcal{V}) \subset T^* G$ is a local diffeomorphism,
and let the quotient maps be given by
$$
\pi_d : G \times G \rightarrow (G \times G)/G \cong G
\quad \text{and} \quad
\pi : T^* G  \rightarrow (T^* G)/G \cong {\mathfrak g}^*.
$$
Let $\ell$ be the reduced Lagrangian on $G$ defined by
$$\mathbb{L} = \ell \circ \pi_d,$$
and let
$$F \ell : \mathcal{U} \subset G \rightarrow {\mathfrak g}^*$$
be the corresponding Legendre transform.
Then the following diagram commutes:
\begin{equation}
\label{diag_com}
\begin{array}{ccc}
\mathcal{V} \subset G \times G & \stackrel{F \mathbb{L}}{\longrightarrow} & T^* G \\
\Big\downarrow\vcenter{%
   \rlap{$\scriptstyle{\rm }\pi_d$}}&
         & \Big\downarrow\vcenter{%
            \rlap{$\scriptstyle{\rm }\pi$}}\\
\mathcal{U} \subset G & \stackrel{F \ell}{\longrightarrow} & {\mathfrak g}^*.
\end{array}
\end{equation}
\end{thm}

\begin{proof}
First, we choose coordinate systems on each space.
Let $(g_k, g_{k+1}) \in G \times G$ and $(g,p) \in T^* G$,
so that the discrete quotient map  (\ref{projection})
is given by $\pi_d : (g_k, g_{k+1}) \mapsto f_{k k+1} = g_k g_{k+1}^{-1}$,
and the continuous quotient map by $\pi : (g,p) \mapsto \mu = R^*_g p$.
Recall that the fiber derivative $F \mathbb{L}$ in these coordinates
has the following form (see, e.g., \cite{WM})
$$
F \mathbb{L} : G \times G \rightarrow T^* G;
\quad (g_k, g_{k+1}) \mapsto (g_k, D_1 \mathbb{L}(g_k, g_{k+1})).
$$

Then the above diagram is given by:
\begin{equation}
\label{diag_coor}
\begin{array}{ccc}
(g_k, g_{k+1}) & \stackrel{F \mathbb{L}}{\longrightarrow} &
(g_k, p_k = \dfrac{\partial \mathbb{L}}{\partial g_k}) \\
\Big\downarrow\vcenter{%
   \rlap{$\scriptstyle{\rm }\pi_d$}}&
         & \Big\downarrow\vcenter{%
            \rlap{$\scriptstyle{\rm }\pi$}}\\
f = R_{g^{-1}_{k+1}} g_k &  &
\mu = R^*_{g_k} p_k,
\end{array}
\end{equation}
where $f$ stands for $f_{k k+1} = g_k g_{k+1}^{-1}$.
To close this diagram and to verify the arrow determined by $F \ell$
compute the derivative of $\mathbb{L}$ using the chain rule
\begin{equation}
   \label{mu_Leg}
   \mu = R^*_{g_k} p_k
   = R^*_{g_k} \frac{\partial (\ell \circ \pi)}{\partial g_k}
   = R^*_{g_k} \left( R^*_{g^{-1}_{k+1}} \frac{\partial \ell}
   {\partial f} \right) = R^*_f \frac{\partial \ell}{\partial f}
   = R^*_f \circ \ell' (f),
\end{equation}
where we have used that, according to the definition of $f$, the
partial derivative $\dfrac{\partial f}{\partial g_k}$ is
given by the linear operator $T R_{g^{-1}_{k=1}}$.
(\ref{mu_Leg}) is precisely the Legendre transformation $F \ell$
for a right invariant system (see the previous subsection).
%and thus it furnishes the proof.
\end{proof}

\begin{cor}
Reconstruction of the discrete Lie-Poisson (DLP) dynamics on
${\mathfrak g}^*$ by $\pi^{-1}$ corresponds to the image of
the discrete Euler-Lagrange (DEL) dynamics on $G \times G$
under the Legendre transformations $F \mathbb{L}$ and
results in an algorithm on $T^*G$ approximating the continuous
flow of the corresponding Hamiltonian system.
\end{cor}

\begin{proof}
The proof follows from the results of the previous subsection,
in particular, diagram (\ref{diag_dyn}) relates  the DLP
dynamics on ${\mathfrak g}^*$ with the DEP dynamics on
$\mathcal{U} \subset G$ which, in turn, is related to the
DEL dynamics on $\mathcal{V} \subset G \times G$ via the
reconstruction (\ref{rec_DEL}).
\end{proof}

An important remark to this corollary which follows from the
results in \cite{KMO}
(see also \cite{KMOW}) is that, in general, to get
a corresponding algorithm on the Hamiltonian side which is
consistent with the corresponding continuous Hamiltonian
system on $T^* G$, one must use the time step $h$-dependent
Legendre transform given by the map
$$
\quad (g_k, g_{k+1}) \mapsto (g_k, -h D_1 \mathbb{L}(g_k, g_{k+1})).
$$
The results of this paper are not effected, however, as we assume
$h$ to be constant and so we would simply add a constant multiplier
to the corresponding symplectic and Poisson structures.
For variable time-stepping algorithms, this remark is crucial and must
be taken into account.

\begin{thm}
The Poisson structure on the Lie group $G$ obtained by  reduction of
the  Lagrange symplectic form $\omega_{\mathbb{L}}$
on $\mathcal{V} \subset G \times G$ via $\pi_d$ coincides with the
Poisson structure on $\mathcal{U} \subset G$ obtained by the
pull-back of the Lie-Poisson structure $\omega_\mu$ on
${\mathfrak g}^*$ by the Legendre transformation $F \ell$
(see diagram (\ref{diag_com}) above).
\end{thm}

\begin{proof}
The proof is based on the commutativity of diagram (\ref{diag_com})
and the $G$ invariance of the unreduced symplectic forms.
Notice that $G$ and ${\mathfrak g}^*$ in (\ref{diag_com}) are
Poisson manifolds, each being foliated by symplectic leaves,
which we denote $\Sigma_f$ and ${\mathcal O}_\mu$ for
$f \in G$ and $\mu \in {\mathfrak g}^*$, respectively.
Denote by $\omega_f$ and $\omega_\mu$ the corresponding
symplectic forms on these leaves. Below we shall prove the
compatibility of these structures under the diagram (\ref{diag_com}).
Repeating this proof leaf by leaf establishes then the equivalence
of the Poisson structures and proves the theorem.

Recall that the Lagrange $2$-form $\omega_{\mathbb{L}}$ on
$\mathcal{V} \subset G\times G$ derived from the variational 
principle coincides
with the pull-back of the canonical $2$-form $\omega_{\text{can}}$
on $T^* G$ (see, e.g., \cite{MPS, WM}).
Recall also that for a right-invariant system, reduction
of $T^* G$ to ${\mathfrak g}^*$ is given by  right
translation to the identity $e \in G$, i.e. any $p \in T^*_g G$ is
mapped to $\mu = R_g^* p \in {\mathfrak g}^* \cong T^*_e G$.
Thus, for any $g \in \pi^{-1}(\mu)$, where $\mu \in {\mathfrak g}^*$,
$$
\pi^{-1} |_{T^*G} = R^*_{g^{-1}} : {\mathfrak g}^* \rightarrow T^*_g G,
$$
so that $(\pi^{-1})^* = (R^*_{g^{-1}})^*$ pulls back
$\omega_{\text{can}}$ to $\omega_\mu$.
Henceforth, $\pi^{-1}$ shall denote the inverse map of $\pi$ restricted
to $T_g G^*$.

Let us write down using the above notations
how the symplectic forms are being mapped
under the transformations in diagram (\ref{diag_com}); we see that
\begin{equation}
\label{diag_symp}
\begin{array}{ccc}
\mathcal{V} \subset G \times G & \stackrel{F \mathbb{L}}{\longrightarrow}
& T^* G \\
\Big\downarrow\vcenter{%
   \rlap{$\scriptstyle{\rm }\pi_d$}}&
         & \Big\downarrow\vcenter{%
            \rlap{$\scriptstyle{\rm }\pi$}}\\
\mathcal{U} \subset G & \stackrel{F \ell}{\longrightarrow} & {\mathfrak g}^*
\end{array}
\qquad \qquad
\begin{array}{ccc}
\omega_{\mathbb{L}} & \stackrel{F \mathbb{L}^*}{\longleftarrow} &
\omega_{\text{can}} \\
%\Big\downarrow\vcenter{%
%  \rlap{$\scriptstyle{\rm }\rho_d$}}
&        & \Big\uparrow\vcenter{%
            \rlap{$\scriptstyle{\rm }(\pi^{-1})^*$}}\\
\omega_f & \stackrel{F \ell^*}{\longleftarrow} & \omega_\mu.
\end{array}
\end{equation}
Then, using the coordinate notations of diagram (\ref{diag_coor}),
for any $f \in G$ and $u, v \in T_f \Sigma_f$,
\begin{equation}
\label{omega1}
\omega_f (f) (u,v) \equiv \omega_\mu (\mu)
(T F \ell (u), T F \ell (v)),
\end{equation}
where $\mu = F \ell (f) \in {\mathfrak g}^*$.
Continuing this equation using diagram (\ref{diag_symp}), we have that
\begin{multline}
\label{omega2}
\omega_f (f) (u,v) = \omega_{\text{can}} ((g_k,p_k))
(T \pi^{-1} \circ T F \ell (u), T \pi^{-1} \circ T F \ell (v)) \\
=  \omega_{\mathbb{L}} ((g_k, g_{k+1}))
(T F \mathbb{L}^{-1} \circ T \pi^{-1} \circ T F \ell (u),
T F \mathbb{L}^{-1} \circ T \pi^{-1} \circ T F \ell (v)),
\end{multline}
where $(g_k, p_k) \in \pi^{-1}(\mu)$ and
$T \pi^{-1}$ denotes $T R^*_{g^{-1}}$.

Using (\ref{diag_com}), it follows that
$$
F \ell \circ \pi_d = \pi \circ F \mathbb{L}
$$
and, hence, for the tangent maps, we have that
$$
T F \ell \circ T \pi_d = T \pi \circ T F \mathbb{L}.
$$
So, if $u, v$ in (\ref{omega1}) are images of some $G$ invariant
vector fields $U, V$ on $\mathcal{V} \subset G \times G$, i.e.
$u = T \pi_d (U), v=T \pi_d (V)$,
then from (\ref{omega2}) it follows that
$$
\omega_f (f) (u,v)=\omega_{\mathbb{L}} ((g_k, g_{k+1})) (U, V),
$$
where $(g_k, g_{k+1})=\pi^{-1}_d (f)$ and
$U, V \in T_{(g_k, g_{k+1})} G \times G$. The last equation precisely
means that $\omega_f$ is the discretely reduced symplectic form,
i.e. the image of $\omega_{\mathbb{L}}$ under the quotient map $\pi_d$.
\end{proof}

Analogous theorems hold for the case of left invariant
systems.

\medskip

\noindent {\bf More General Configuration Spaces.}
Similar ideas carry over to the integration of systems defined
on a general configuration space $M$ with some symmetry group $G$.
In this case, the reduced discrete space $(M \times M)/G$ inherits a
Poisson structure from the one defined on $M \times M$
(analogously to (\ref{p2})).
Its symplectic leaves again become dynamically invariant manifolds
for structure-preserving integrators and can be viewed as
images of the symplectic leaves of the reduced Poisson manifold
$T^*M/G$ under appropriately defined ``Legendre transformations''.
This is a topic of ongoing research that builds on recent progress in
Lagrangian reduction theory; see \cite{MRS 99}.

%%%%%%%%%%%%%%%%%%%%%%%%%%%%%%%%%%%%%%%%%%%%%%%%%%%%%%%%%
\subsection{Poisson structures of the rigid body}
%%%%%%%%%%%%%%%%%%%%%%%%%%%%%%%%%%%%%%%%%%%%%%%%%%%%%%%%%

As an example of applications of the above ideas, we consider
the dynamics of the rigid body and its associated reduction
and discretization (see, e.g. \cite{MaRa1999, MoV, LS, MPeS}
for more details).
\medskip

\noindent {\bf The Basic Set Up.} The configuration space of
the system is
$\operatorname{SO}(3)$. The corresponding Lagrangian is
determined by a symmetric positive definite operator
$J : {\mathfrak so}(3) \rightarrow {\mathfrak so}(3)$, defined by
$J (\xi) = \Lambda \xi + \xi \Lambda$, where $\xi \in {\mathfrak so}(3)$
and $\Lambda$ is a diagonal matrix satisfying $\Lambda_i + \Lambda_j>0$
for all $i \ne j$. The left invariant metric on
$\operatorname{SO}(3)$ is obtained by left translating the bilinear
form at $e \in \operatorname{SO}(3)$ given by
$$
( \xi, \xi ) = \frac{1}{4} \operatorname{Tr} \left( \xi^T J (\xi)
\right).
$$

The operator $J$, viewed as a mapping ${\mathcal J} : {\mathfrak so}(3)
\rightarrow {\mathfrak so}(3)^\ast$, has the usual interpretation of
the inertia tensor, and the
$\Lambda_i$ correspond to the sums of certain principal moments of
inertia.

The rigid body Lagrangian is the kinetic energy of the
system
\begin{equation}
\label{rb_Lag}
L (g, \dot{g}) = \frac{1}{4} \langle  g^{-1} \dot{g},
{\mathcal J} ( g^{-1}\dot{g}) \rangle = \frac{1}{4} \langle \xi,
{\mathcal J} \xi) \rangle = l (\xi),
\end{equation}
where $\xi =  g^{-1}\dot{g} \in {\mathfrak so}(3)$ and
$\langle \cdot , \cdot \rangle$
is the pairing between the Lie group and its dual.
\medskip

\noindent {\bf Poisson Structures and Casimir Functions.}
The Lie algebra dual ${\mathfrak so}(3)^\ast$  has a well-known
Lie-Poisson structure with a Casimir
$C_{\mathfrak{so}(3)^\ast} (\mu) = \operatorname{Tr} (\mu^2)$,
where $\mu \in {\mathfrak so}(3)^\ast$.
Upon identification with $\mathbb{R}^3$, its generic
symplectic leaves become concentric spheres with Kirilov-Kostant
symplectic form being proportional to the area form.
If $y$ denotes coordinates on
$\mathbb{R}^3 \cong {\mathfrak so}(3)^\ast$, then the above
Casimir function is given by
$C_{{\mathfrak so}(3)^\ast} (y) = \|y\|^2$.

Following Section 5 of \cite{MPeS} on discrete Euler-Poincar\'{e}
reduction, we obtain the reduced form of the Moser-Veselov
Lagrangian on the group SO$(3)$ given by
$$ \ell (f) = \operatorname{Tr} (f \Lambda), $$
where $f \in \operatorname{SO}(3)$ and
$\operatorname{SO}(3)$ is embedded into the linear space $\mathfrak{gl}(3)$.
Then, the Legendre transform $F \ell$ takes the form
$$
F \ell (f) = L_f^* \circ d \ell (f) =
\operatorname{skew} (f \Lambda) = f \Lambda - \Lambda f^T \
:  \operatorname{SO}(3) \rightarrow {\mathfrak so}(3)^\ast,
$$
where the constraint that $f$ be in $\operatorname{SO}(3)$ has been
enforced.
The pull-back of $C_{{\mathfrak so}(3)^\ast}$ under $F \ell^*$
defines a Casimir function {\it on the group}, which up to a
constant term and a sign, is given by
\begin{equation}
\label{SO_Cas}
C_{\operatorname{SO}(3)} (f) = \operatorname{Tr} (f \Lambda f \Lambda)
\qquad f \in \operatorname{SO}(3).
\end{equation}
Its symplectic leaves constitute the invariant manifolds
of the reduced discrete dynamics corresponding to the
Lagrangian (\ref{rb_Lag}).

%\begin{rem}
Analogously, one can define a Poisson structure on the Lie algebra
${\mathfrak so}(3)$ using the duality between Lie-Poisson and
Euler-Poincar\'{e} reduced systems on ${\mathfrak so}(3)^*$ and
${\mathfrak so}(3)$, respectively. The Lagrangian (\ref{rb_Lag})
defines the Legendre transformations $F l$ from ${\mathfrak so}(3)$
to ${\mathfrak so}(3)^*$ given by
$\mu = \dfrac{\partial l}{\partial \xi} = {\mathcal J} (\xi)$.
Then, the pull-back by $F l^*$ defines a Casimir function on
${\mathfrak so}(3)$:
$$
C_{{\mathfrak so}(3)} (\xi) = F l^* \cdot C_{{\mathfrak so}(3)^*}(\xi)
= \langle \langle {\mathcal J}(\xi), {\mathcal J}(\xi) \rangle \rangle,
$$
where the metric on the dual is induced by the one on the algebra,
i.e. by the symmetric positive definite operator $J$.
If $x$ denotes coordinates on $\mathbb{R}^3 \cong {\mathfrak so}(3)$,
then the above Casimir function is given by
$C_{{\mathfrak so}(3)} (x) = \|{\mathcal J} (x)\|^2$.
Thus, the corresponding symplectic leaves are ellipsoids of
${\mathcal J}^2$. They {\it do not} coincide with adjoint
orbits, which are spheres in $\mathbb{R}^3$.
The dynamic orbits are obtained by intersecting these ellipsoids,
determined by ${\mathcal J}^2$,
with the energy ellipsoids, determined by ${\mathcal J}$.
%\end{rem}

%%%%%%%%%%%%%%%%%%%%%%%%%%%%%%%%%%%%%%%%%%%%%%%%%%%%%%%%%
\section*{Acknowledgments}
%%%%%%%%%%%%%%%%%%%%%%%%%%%%%%%%%%%%%%%%%%%%%%%%%%%%%%%%%
The authors would like to thank Alan Weinstein for pointing
out the connections with the general theory of dynamics on
groupoids and algebroids and Sameer Jalnapurkar for comments
and discussions. SP and SS would also like to thank
the Center for Nonlinear Science in Los Alamos for providing
a valuable setting for part of this work. SS was partially supported by
NSF-KDI grant ATM-98-73133, and JEM and SP were partially supported
by NSF-KDI grant ATM-98-73133 and the Air Force Office of Scientific
Research. 

%%%%%%%%%%%%%%%%%%%%%%%%%%%%%%%%%%%%%%%%%%%%%%%%%%%%%%%%%

\end{document}